# Identities for the String-Product of Strings on Geometric Sequences

**Gbeddy Marcel Selase**†
**(*Wodem*)**

Department of Mathematics
Kwame Nkrumah University of Science and Technology, Kumasi, Ghana
† msgbeddy@hotmail.com

**ABSTRACT**

Geometric sequences are found documented as early as 300BC in the text, Book IX of Elements written by Euclid of Alexandria. In this paper a new principle for identities involving the product of any $k-$ number of terms of a geometric sequence is presented. The results of which are independent from a specific positional term or the common ratio.

## 1- INTRODUCTION

Let $f: X \to Y$ be a one-to-one and onto function such that $f(n) = a_n$ where $X = \{n \in \mathbb{N}, 1 \leq n \leq l\}$ for arbitrarily large $l$. If $a_1, r \in \mathbb{R}:\ a_n = a_{n-1} \cdot r, \forall n \geq 2$ then

$$a_1, a_2, a_3, \ldots, a_l \qquad \text{is a geometric sequence}$$

All closed formula related to a geometric sequence rely explicitly on at least one positional term, $a_m$ and the common ratio $r$. For general notation and results related to geometric sequences see [1].

Let $tuples[Y, t]$ be the list of all $t-$tuples on $Y = \{a_1, a_2, a_3, \ldots, a_l\}$. Denote the different $t-$tuples on $Y$ as the strings $\boldsymbol{S}_i^t$, $1 \leq i \leq \binom{l}{t}$. The string-product of $\boldsymbol{S}_k^t$ is $\prod_{j \in \{f^{-1}: \boldsymbol{S}_k^t \to X\}} a_j$, the product of all elements in $\boldsymbol{S}_k^t$. The sum of subscripts of elements in $\boldsymbol{S}_k^t$ is also the sum of objects in $X$ that have an image in $\boldsymbol{S}_k^t$ under $f$, that is $\sum_{j \in \{f^{-1}: \boldsymbol{S}_k^t \to X\}} j$.

## 2- String-Product of Strings

**Lemma 2.1** For a general geometric sequence $a_1, a_2, a_3, \dots, a_l$ the $n$th term $a_n$ given as $a_n = a_1 \cdot r^{n-1}$.

*Proof.* $a_n = a_{n-1} \cdot r = a_{n-2} \cdot r^2 = a_{n-3} \cdot r^3 = \cdots = a_{n-k} \cdot r^k$ for $n > k \in \mathbb{N}$. If $k = n - 1$, we have $a_n = a_{n-(n-1)} \cdot r^{n-1} = a_1 \cdot r^{n-1}$

**Corollary 2.2** If $\boldsymbol{S}_u^t$ and $\boldsymbol{S}_v^t$ are $t$ –tuples on $Y$, such that $\prod_{j\in\{f^{-1}:S_u^t\to X\}} a_j = \prod_{k\in\{f^{-1}:S_v^t\to X\}} a_k$, then $\sum_{j\in\{f^{-1}:S_u^t\to X\}} j = \sum_{k\in\{f^{-1}:S_v^t\to X\}} k$

*Proof.* The result is a direct consequence of Lemma 2.1 since, if $\prod_{j\in\{f^{-1}:S_u^t\to X\}} a_j = \prod_{k\in\{f^{-1}:S_v^t\to X\}} a_k$, then $\prod_{j\in\{f^{-1}:S_u^t\to X\}}[a_1 \cdot r^{j-1}] = \prod_{k\in\{f^{-1}:S_v^t\to X\}}[a_1 \cdot r^{k-1}]$

This gives $\prod_{j\in\{f^{-1}:S_u^t\to X\}} a_1 \cdot \prod_{j\in\{f^{-1}:S_u^t\to X\}} r^j \cdot \prod_{j\in\{f^{-1}:S_u^t\to X\}} r^{-1} = \prod_{k\in\{f^{-1}:S_v^t\to X\}} a_1 \cdot \prod_{k\in\{f^{-1}:S_v^t\to X\}} r^k \cdot \prod_{k\in\{f^{-1}:S_v^t\to X\}} r^{-1}$

Since both $\boldsymbol{S}_u^t$ and $\boldsymbol{S}_v^t$ are $t$ –tuples, $\prod_{j\in\{f^{-1}:S_u^t\to X\}} a_1 = \prod_{k\in\{f^{-1}:S_v^t\to X\}} a_1 = (a_1)^t$ and $\prod_{j\in\{f^{-1}:S_u^t\to X\}} r^{-1} = \prod_{k\in\{f^{-1}:S_v^t\to X\}} r^{-1} = r^{-t}$. The equation therefore reduces to

$\prod_{j\in\{f^{-1}:S_u^t\to X\}} r^j = \prod_{k\in\{f^{-1}:S_v^t\to X\}} r^k$ hence the result $r^{\sum_{j\in\{f^{-1}:S_u^t\to X\}} j} = r^{\sum_{k\in\{f^{-1}:S_v^t\to X\}} k}$

$\sum_{j\in\{f^{-1}:S_u^t\to X\}} j = \sum_{k\in\{f^{-1}:S_v^t\to X\}} k$

The result leads to the following theorem.

**Theorem 2.3** For any geometric sequence $a_1, a_2, a_3, \dots, a_l$, where $a_n$ is the $nth$ term, the product of the same number of terms is the same provided the sum of the subscripts are the same.

*Proof.* The result is a direct consequence of the method of proof for Corollary 2.2 and Lemma 2.1.

**Example 2.3.1** If we take two terms like $a_3$ and $a_4$ we can write $a_4 \cdot a_3 = a_5 \cdot a_2 = a_6 \cdot a_1$ since the sum of subscripts in each pair is 7. Generally for two terms $a_i \cdot a_j = a_{i-n} \cdot a_{j+n} \; \forall \; i < n \; where \; i, j, n \in \mathbb{N}$. If we take three terms like $a_2, a_3$ and $a_7$, we can write $a_2 \cdot a_3 \cdot a_7 = a_4 \cdot a_2 \cdot a_6 = a_5 \cdot a_3 \cdot a_4 = a_8 \cdot a_3 \cdot a_1 = \cdots$, as long as the sum of the subscripts of the three terms is 12.

**Lemma 2.4** If $\sum_{j\in\{f^{-1}:S_u^t\to X\}} j$ can be written as $t_1b_1 + t_2b_2 + \cdots + t_nb_n$, $\forall b_i \in \mathbb{N}, 1 \leq b_i \leq l$ where $t_1 + t_2 + \ldots + t_n = t$, then $\sum_{j\in\{f^{-1}:S_k^t\to X\}} a_j = \left(a_{b_1}\right)^{t_1} \cdot \left(a_{b_2}\right)^{t_2} \cdot \ldots \cdot \left(a_{b_n}\right)^{t_n}$.

*Proof*. By the associative law of multiplication and immediate induction thereafter, it suffices to prove the result for a 2-tuple $\{a_i, a_j\}$ on $Y$.

Assume $i + j = 2k$. Thus $a_i \cdot a_j = a_1 \cdot r^{n-1} \cdot a_1 \cdot r^{n-1} = (a_1)^2 \cdot r^{i+j-2} = (a_1)^2 \cdot r^{2k-2} = (a_1 \cdot r^{k-1})^2 = (a_k)^2$. Therefore the result in general

**Corollary 2.5** If in Lemma 2.4 $\frac{t_1}{n_1} + \frac{t_2}{n_2} = t$ and $\frac{t_1}{n_1} \cdot i + \frac{t_2}{n_2} \cdot j = t \cdot k$ then for $\{a_i, a_j\}$ $(a_i)^{\frac{t_1}{n_1}} \cdot (a_j)^{\frac{t_2}{n_2}} = (a_k)^t$

*Proof.* The result is a direct consequence of the method of proof of Lemma 2.4 together with the laws for exponents which are not restricted to natural numbers. And it remains valid for $t_1, t_2 \in \mathbb{R}$.

**Example 2.5.1** For any geometric sequence consider the terms $a_5, a_2$ and $a_4$.
Since $1 + \frac{1}{2} = \frac{3}{2}$ and $1(5) + \frac{1}{2}(2) = \frac{3}{2}(4)$ we can write $a_5 \cdot (a_2)^{\frac{1}{2}} = (a_4)^{\frac{3}{2}}$

Also since $6\pi + 6 = (5\pi + 2) + (\pi + 4)$ and $6\pi(3) + 6(6) = (5\pi + 2)(2) + (\pi + 4)(8)$ we can write $(a_3)^{6\pi} \cdot (a_6)^6 = (a_2)^{5\pi+2} \cdot (a_8)^{\pi+4}$

## 3- CONCLUSION

It was found that applying the results to geometric sequences in general gets cumbersome with increase in the number terms used. This observation leaves room for further research to find results which will ease the applications for large number of terms.

***Acknowledgement***: The author gratefully acknowledge the critical and creative suggestions and guidelines by Dr Johan Kok, who is an independent researcher from South Africa and is a specialist research advisor for the Centre of Studies in Discrete Mathematics, Thrissur, Kerala, India.

## 4- REFERNCES